%% Amstex file, this is the version of July 2013.
%% August 8, 2013. Accepted for publication in
%% "Journal of Analysis & Number Theory". 18.09.2013.

\magnification1200
\input amstex.tex
\documentstyle{amsppt}

\hsize=12.5cm
\vsize=18cm
\hoffset=1cm
\voffset=2cm

\footline={\hss{\vbox to 2cm{\vfil\hbox{\rm\folio}}}\hss}
\nopagenumbers
\def\DJ{\leavevmode\setbox0=\hbox{D}\kern0pt\rlap
{\kern.04em\raise.188\ht0\hbox{-}}D}

\def\txt#1{{\textstyle{#1}}}
\baselineskip=13pt
\def\hf{{\textstyle{1\over2}}}
\def\a{\alpha}\def\b{\beta}
\def\d{{\,\roman d}}
\def\e{\varepsilon}

\def\G{\Gamma}

\def\={\;=\;}

\def\zt{\zeta(\hf+it)}

\def\D{\Delta}

\def\R{\Re{\roman e}\,} 
\def\z{\zeta}

\def\hf{{\textstyle{1\over2}}}
\def\txt#1{{\textstyle{#1}}}

\def\le{\leqslant} \def\ge{\geqslant}
%%%%%%%%%%% Fonts macros %%%%%%%%%%%%
\font\tenmsb=msbm10
\font\sevenmsb=msbm7
\font\fivemsb=msbm5
\newfam\msbfam
\textfont\msbfam=\tenmsb
\scriptfont\msbfam=\sevenmsb
\scriptscriptfont\msbfam=\fivemsb
\def\Bbb#1{{\fam\msbfam #1}}

\def \NN {\Bbb N}

\font\ff=cmr8
\def\txt#1{{\textstyle{#1}}}
\baselineskip=13pt

\font\teneufm=eufm10
\font\seveneufm=eufm7
\font\fiveeufm=eufm5
\newfam\eufmfam
\textfont\eufmfam=\teneufm
\scriptfont\eufmfam=\seveneufm
\scriptscriptfont\eufmfam=\fiveeufm
\def\mathfrak#1{{\fam\eufmfam\relax#1}}

\font\tenmsb=msbm10
\font\sevenmsb=msbm7
\font\fivemsb=msbm5
\newfam\msbfam
     \textfont\msbfam=\tenmsb
      \scriptfont\msbfam=\sevenmsb
      \scriptscriptfont\msbfam=\fivemsb
\def\Bbb#1{{\fam\msbfam #1}}

\def \NN {\Bbb N}

  \def\rightheadline{{\hfil{\ff
  The divisor function and divisor problem}\hfil\tenrm\folio}}

  \def\leftheadline{{\tenrm\folio\hfil{\ff
   Aleksandar Ivi\'c }\hfil}}
  \def\emptyheadline{\hfil}
  \headline{\ifnum\pageno=1 \emptyheadline\else
  \ifodd\pageno \rightheadline \else \leftheadline\fi\fi}

\topmatter
\title
THE DIVISOR FUNCTION AND DIVISOR PROBLEM
\endtitle
\author   Aleksandar Ivi\'c  \endauthor

\bigskip
%\centerline{\bb }
%\bigskip
\address
Aleksandar Ivi\'c, Katedra Matematike RGF-a
Universiteta u Beogradu, \DJ u\v sina 7, 11000 Beograd, Serbia
\endaddress

\email {\tt
ivic\@rgf.bg.ac.rs,  aivic\@matf.bg.ac.rs} \endemail
\dedicatory
Extended version of the invited talk given in Budapest at the ERDOS100 Conference, July 1-5, 2013
\enddedicatory
\abstract
{The purpose of this text is twofold. First we discuss some problems
involving Paul Erd\H os (1913-1996), whose centenary of birth is this year. In the second part
some recent results on divisor problems are discussed, and their connection
with the powers moments of $|\zt|$ is pointed out. }
\endabstract
\endtopmatter

\document

\head
1. Introduction
\endhead
The classical number of divisors function of a positive integer $n$ is
$$
d(n) := \sum_{\delta|n}1.
$$
We have
$d(mn) = d(m)d(n)
$
whenever $(m,n)=1$, so that $d(n)$ is a
multiplicative arithmetic function. Further $d(p^\a)=\a+1$
for $\a\in\NN$, where  $p, p_j$ denote generic primes and $\NN$ is the set of natural
numbers. Therefore, if
$n = \prod_{j=1}^r p_j^{\a_j}$ is the canonical decomposition of $n$ into
prime powers, then
$$
d(n) = (\a_1+1)(\a_2+1)\cdots (\a_r+1).
$$
In general
$$
\z^k(s) = \sum_{n=1}^\infty d_k(n)n^{-s}\qquad(k\in\NN,\;\R s>1),
$$
where the (general) divisor function $d_k(n)$ represents the number
of ways $n$ can be written as a product of
$k$ factors, so that in particular $d_1(n)\equiv 1$ and $d(n) \equiv d_2(n)$.
The Riemann zeta-function is
$$
\z(s) := \sum_{n=1}^\infty n^{-s} = \prod_p{(1-p^{-s})}^{-1}\qquad (\R s > 1),
$$
otherwise it is defined by analytic continuation. It is the simplest and most important
example of the so-called class of $L$-functions $\sum_{n=1}^\infty f(n)n^{-s}$ satisfying certain natural
properties. See e.g., J. Kaczorowski [Kac] for a survey of the {\it Selberg class} of $L$-functions
generalizing $\z(s)$.
\medskip
The function $d_k(n)$ is a also multiplicative
function of $n$, meaning that $d_k(mn) = d_k(m)d_k(n)$ if $m$ and $n$
($\in\NN$) are coprime, and
$$
d_k(p^\a) = (-1)^\a{\binom{-k}{\a}} = {k(k+1)\cdots (k+\a-1)\over \a!}
$$
for primes $p$ and $\a\in\NN$.

\medskip
\head
2. Iterations of $d(n)$
\endhead
From the wealth of problems involving the divisor function $d(m)$ we shall concentrate
on some problems connected with the work of Paul Erd\H os (1913-1996), one of the greatest
mathematicians of the XXth century. We begin with the iterations
of $d(n)$. Thus let, for $k\in\NN$ fixed,
$$
d^{(1)}(n) : = d(n),\; d^{(k)}(n) := d\left(d^{(k-1)}(n)\right)\quad(k>1)
$$
be the $k$-th iteration of $d(n)$. Already $d^{(2)}(n)$ {\bf is not}
multiplicative! This fact makes the problems involving $d^{(k)}(n)$ and
iterates of other multiplicative functions quite difficult.

\medskip
The great Indian mathematician S. Ramanujan (1887-1920) [Ram] proved in 1915 that
(for the connection between Erd\H os and Ramanujan see [Erd3])
$$
d^{(2)}(n) \;>\; 4^{\sqrt{2\log n}/\log\log n}
%\qquad(n\geq n_0)
$$
for infinitely many $n$.
This lower bound follows if one considers ($p_j$ is the $j$-th prime)
$$
N = 2^1\cdot 3^2\cdot 5^4\cdot \ldots \cdot p_k^{p_k-1}\leqno(2.1)
$$
and lets $k\to\infty$. Namely
$$
d(N) = 2\cdot3\cdot5\cdot \ldots \cdot p_k, \quad d^{(2)}(N) = 2^k,
$$
and one easily bounds $k$ from (2.1) by the prime number theorem. Ramanujan's
paper [Ram] contains many other results on the divisor function $d(n)$.

\medskip
Important results on the order of $d^{(k)}(n)$ were obtained in 1967 by P.
 Erd\H os and I. K\'atai [ErKa].
Let $\ell_k$ denote the $k$-th Fibonacci number:
$$
\ell_{-1}=0,\; \ell_0=1,\; \ell_k = \ell_{k-1}+\ell_{k-2}\quad(k\geqslant1).
$$
Then the result of P. Erd\H os and I. K\'atai says that
$$
d^{(k)}(n) \;<\; \exp\left((\log n)^{1/\ell_k+\e}\right)\leqno(2.2)
$$
for fixed $k$ and $n\geqslant n_0(\e,k)$, and that for every $\e>0$
$$
d^{(k)}(n) \;>\; \exp\left((\log n)^{1/\ell_k-\e}\right)\leqno(2.3)
$$
for infinitely many $n$. Here and later
 $\e$ denotes arbitrarily small positive constants, not
 necessarily the same ones at each occurrence.
The lower bound in (2.3) follows for $n = N_j$, where inductively
$N_1 = 2\cdot 3\cdot\ldots\cdot p_r,$
and if
$$
N_j = \prod_{i=1}^{S_j}p_i^{r_i},
$$
say, then
$$
\eqalign{
N_{j+1} &= (p_1\cdot \ldots \cdot p_{r_1})^{p_i-1}(p_{r_i+1}\cdot \ldots \cdot
 p_{r_1+r_2})^{p_2-1}\cr
&\cdots (p_{r_1+\cdots r_{S_j-1}+1}\cdot \ldots \cdot p_{r_1+\cdots r_{S_j}})^{p_{S_j-1}}.\cr}
$$
Then one has $d^{(k)}(N_k) = 2^r$, and the proof reduces to finding
the lower bound for $r$. The proof of the upper bound in (2.2) is more involved.

\medskip

Improvements of (2.2) and (2.3) in the general case have been obtained 
by A. Smati [Sm2]. Extensive work has been done in the case when $k=2$.
Thus in [ErIv1] P. Erd\H os and A. Ivi\'c
proved, for $n\geqslant n_0$ and suitable $ C >0$,
$$
d^{(2)}(n) < \exp\left(C\left(\frac{\log n\log\log n}{\log\log\log n}\right)^{1/2}\right).
\leqno(2.4)
$$

This follows from
$$
\log d(n) = \sum_{i=1}^r\log(\a_i+1) \ll r\log\log r = \omega(n)\log\log \omega(n),
$$
and the bound ($\omega(n) = \sum_{p|n}1$ is the number of distinct prime factors of $n >1$, $\omega(1)=0$)
$$
\omega(d(n)) \ll \left(\frac{\log n\log\log n}{\log\log\log n}\right)^{1/2}.
$$
For some other work of P. Erd\H os and A. Ivi\'c on $d(n)$ see [ErIv2] and [EGIP], the
last paper bring a joint work with S.W. Graham and C. Pomerance.

\medskip
A. Smati [Sm1], [Sm2] improved the  upper bound
in (2.4) to

$$
d^{(2)}(n) < \exp\bigl(C\sqrt{\log n}\,\bigr)\qquad(C>0, n \geqslant n_0),\leqno(2.5)
$$
 which turned out to be only by a factor of $\log\log n$
(in the exponent) smaller
than the true upper bound.
Namely, in 2011 Y. Buttkewitz, C. Elsholtz, K. Ford and J.-C. Schlage-Puchta [BEFSP] practically
settled the problem of the maximal order of $d^{(2)}(n)$  by proving that
$$
\max_{n\leqslant x}\log \,d^{(2)}(n) \;=\; \frac{\sqrt{\log x}}{\log\log x}\left(D+
O\left(\frac{\log\log\log x}{\log\log x}\right)\right),
$$
where $D= 2.7958\ldots\;$ is an explicit constant.
Note: We use throughout the paper the notation ($C$ denotes generic postive constants)
$$
f(x) \ll g(x) \Longleftrightarrow f(x) = O(g(x)) \Longleftrightarrow
|f(x)| \leqslant Cg(x)\quad(x\geqslant x_0).
$$

\medskip
R. Bellman and H.N. Shapiro [BeSh] conjectured that, for fixed $k\ge 1$,
$$
\sum_{n\le x}d^{(k)}(n) = (1+o(1))c_kx\log_kx\qquad(x\to\infty),
\leqno(2.6)
$$
where $\log_k$ is the $k$ times iterated natural logarithm. For $k=1$ this is trivial,
but for $k>1$ it is a difficult problem. P. Erd\H os [Erd2] and I. K\'atai [Ka1]
obtained (2.6) for $k=2$,
while I. K\'atai [Ka2] proved it for $k=3$. Finally
 Erd\H os and K\'atai [ErKa2] proved it for $k=4$, where the matter seems to stand at present.

\medskip
Finally we mention a problem related to the iteration of $d(n)$.
In 1992 the author [Iv2]  conjectured that
$$
\sum_{n\leqslant x}d\bigl(n+d(n)\bigr)= Bx\log x + O(x)\quad(B>0). \leqno(2.7)
$$
I. K\'atai [Ka3] obtained this formula with the error term
$O(x\log x/\log\log x)$. He indicated that a formula analogous to (2.7) holds
also for the summatory function of $d(n+f(n))$, where $f(n) = \omega(n), \Omega(n) = \sum_{p^\a||n}\a,
d_k(n)$, for example.

\medskip
\head
3. P. Erd\H os's work on $d(n)$ in short intervals
\endhead

From the rich legacy of P. Erd\H os concerning results and problems involving $d(n)$
we single out  his classical paper  [Erd1] (for some of his other papers
involving $d(n)$, see  [ENS] (with J.-L. Nicolas and A. S\'ark\"ozy)
and [ErHa] and [ErMi], the last two written jointly with R.R. Hall
and L. Mirsky, respectively). He begins in [Erd1] (we keep his German
original):

\medskip

$d(n)$ {\it sei der Anzahl der Teiler von $n$. Folgende asymptotische Formel ist wohl-\break bekannt:
$$
\sum_{n=1}^x d(n) = x\log x + (2C-1)x + O(x^\a), \quad\a = 15/46\leqno(3.1)
$$
($C$ ist die Eulersche Konstante)}.

\medskip
Note that the function in the $O$-term is commonly denoted by $\D(x)$, thus
$$
\D(x) := \sum_{n\leqslant x}d(n) - x(\log x + 2C -1).\leqno(3.2)
$$
The constant 15/46 = 0.32608$\ldots$, due to H.-E. Richert [Ric] (1952), can be replaced by M.N. Huxley's
(2003)[Hux] value 131/416 = 0.31493$\ldots\;$.

\medskip
Erd\H os's theorem is as follows: {\it Es sei $h(x)$ eine beliebige wachsende Funktion,
die mit $x$ gegen $\infty$ strebt. Es sei
$$
f(x) > (\log x)^{2\log2-1}\exp\Bigl(h(x)\sqrt{\log\log x}\,\Bigr).
$$
Dann gilt f\"ur fast alle $x$
$$
\sum_{n\leqslant f(x)}d(x+n) = (1+o(1))f(x)\log x\qquad(x\to\infty).\leqno(3.3)
$$

\medskip
Diese Formel l\"asst sich nicht weiter versch\"arfen. Ist n\"amlich
$$
f(x) =  (\log x)^{2\log2-1}\exp\Bigl(c\sqrt{\log\log x}\,\Bigr)\qquad(c>0),
$$
so gilt (3.3) nicht mehr f\"ur fast alle $x$.}

\medskip
It is commonly conjectured that the error term in (3.1) is $O_\e(x^{1/4+\e})$,
while it is known long ago that it is $\Omega(x^{1/4})$ $(f(x) = \Omega(g(x))$ as
$x\to\infty$ means that $f(x) = o(g(x))$ does not hold). The conjecture  on the
upper bound is one of the most difficult problems in analytic number theory, as
it does not appear to follow from the Lindel\"of hypothesis (LH, $\zt \ll_\e |t|^\e$) or
from the Riemann hypothesis (RH, all complex zeros of $\z(s)$ satisfy $\R s = \hf$).
Both the LH and RH are unsettled to this day, and it is known that the RH
implies the LH (see e.g., Chapter 1 of [Iv1]).

\medskip
Let $D_k^+(n)=\max_{0\leq h<k}d(n+h)$, where $d(n)$ is the number of divisors of $n$. 
P. Erd\H os and R.R. Hall [ErHa]  showed that for fixed $k$, 
$$
\sum_{n\leq x}D_k^+(n) = k(1+o(1))x\log x\qquad(x\rightarrow\infty).\leqno(3.4)
$$
 R.R. Hall conjectured that (3.4) is true so long as 
$k\leq(\log x)^\alpha$ with $\alpha<\log 4-1$. He showed that (3.4) fails if $k\geq(\log x)^\alpha$ with $\alpha>\log 4-1$. 
R.R. Hall [Hal] showed later that in fact (3.4) is true for 
$k\leq(\log x)^{\log 4-1}\exp\{-\xi(x)\root{}\of{\log\text{}\log x}\,\}$ where 
$\xi(x)\rightarrow\infty$. Furthermore $\sum_{n\leq x}D_k^+(n)=o(kx\log x)$ 
if $$
k>(\log x)^{\log 4-1}\exp\Bigl\{\xi(x)\root{}\of{\log\text{}\log x}\,\Bigr\}.
$$ 
Further results on this and on related topics were obtained by R.R. Hall and G. Tenenbaum [HaTe].
\medskip
\head
4. The additive divisor problem
\endhead
\medskip
We turn now to modern developments involving the divisor function in short intervals.
The importance of these results is that they have applications to power moments of $|\zt|$, which is
one of central topics in the theory of the Riemann zeta-function.
The author [Iv4]  proved in 1997 the folllowing

\medskip
THEOREM 1. {\it
 For a fixed integer $k \ge 3$ and any fixed $\e > 0$,
we have
$$
\int_0^T|\zeta(\hf+it)|^{2k}\d t \;\ll_{k,\e}\; T^{1+\e}\Bigl(1 +
\sup_{T^{1+\e}<M\ll T^{k/2}}\,
G_k(M;T)M^{-1}\Bigr), \leqno(4.1)
$$
if, for $M < M' \le 2M, T^{1+\e} \le M \ll T^{k/2}$,
$$
G_k(M;T)
:= \sup_{\scriptstyle M\le x \le M'\atop\scriptstyle 1\le t \le M^{1+\e}/T}
\Bigl|\sum_{h\le t}\,\D _k(x,h)\Bigr|.
$$
}
\medskip
The bound (4.1) provides a direct
link between upper bounds for the $2k$-th
moment of $|\zt|$
and sums of $\D _k(x,h)$ over the shift parameter $h$, showing also the limitations
of the method, where $\D _k(x,h)$ denotes the error term in the
asymptotic formula for the sum $\sum_{n\leqslant x}d_k(n)d_k(n+h)$.
Of course the problem greatly increases in complexity as $k$ increases,
and this is one of the reasons why in [Iv3] only the case $k =3$ was
considered. The case $k = 2$ was not treated at all, since for the fourth
moment of $|\zt|$ we have an asymptotic formula with precise results for
the corresponding error term (see e.g., Chapter 5 of [Iv1] and the paper
of Y. Motohashi and the author [IvMo1]).

As for the function $\Delta_k(x,h)$,  one writes
$$
\sum_{n \le x} d_k(n) d_k(n+h) = x\,P_{2k-2}(\log x;h)+\Delta_k(x,h),
$$
where it is assumed that $k\ge2$ is a fixed integer,
and $P_{2k-2}(\log x;h)$ is a suitable polynomial of degree $2k-2$ in $\log x$,
whose coefficients depend on $k$ and $h$,
while $\Delta_k(x,h)$ is supposed to be the error term.
This means that we should have
$$
\Delta_k(x,h)\; =\; o(x) \quad{\roman as}\quad x \to \infty,
$$
but unfortunately this is not yet known to hold for any $k \ge 3$, even for fixed $h$,
while for $k=2$ there are many results. This is the so-called {\it binary additive divisor
problem} in the case when $k=2$, and the {\it general additive divisor problem}
when $k>2$. The binary additive divisor problem consists of the evaluation of the sum
$$
D(N;f) := \sum_{n\le N}d(n)d(n+f),\leqno(4.2)
$$
where $f$ is a natural number, not necessarily fixed. One can write
$$
D(N;f) = M(N;f) + E(N;f),\leqno(4.3)
$$
where $M(N;f)$ and $E(N;f)$ are to be considered the ``main term'' and the ``error term'',
respectively, in the asymptotic formula for the sum $D(N;f)$ in (4.2) as $N\to\infty$.
Already A.E. Ingham [Ing] showed that the main term in (4.3) has the form
$$
M(N;f) = \Bigl\{c_1(f)\log^2N + c_2(f)\log N + c_3(f)\Bigr\}N,
$$
where the coefficients $c_j(f)$ (which depend on $f$) can be written down explicitly.
For some modern results on $E(N;f)$ we refer the reader to [Mot], [IvMo2] and {Meu].

\medskip
It seems reasonable to expect that for the quantity $G_k$ in (4.1)
we shall have a bound of the form
$$
G_k \,\ll_{k,\e} \,T^{a_k+\e}M^{b_k+\e} \qquad (a_k \ge 0, b_k \ge 1)
\leqno(4.4)
$$
with suitable constants $a_k,b_k$. Hence assuming (4.4) it follows that
$$
G_kM^{-1} \ll T^{a_k+\e}M^{(b_k-1+\e)} \ll T^{a_k+{k\over2}(b_k-1)+\e}
$$
for $M \ll T^{k/2}$. Therefore from Theorem 1 we obtain

\medskip
{\bf Corollary 1}. {\it If (4.4) holds, then for
a fixed integer $k \ge 3$ we have}
$$
\int_0^T|\zt|^{2k}\,\d t\;\ll_{k,\e}\;T^{1+\e}\left(1 +
T^{a_k+{k\over2}(b_k-1)}\right).\leqno(4.5)
$$

\medskip
The case $k=3$ was investigated by the author in [Iv4].
Therein it was conjectured  that
$$
\sum_{h\le H}\D _3(x,h) \ll_\e Hx^{{2\over3}+\e}\qquad(1 \le H \le x^{{
1\over3}+\delta}) \leqno(4.6)
$$
for some $\delta > 0$. If $k = 3$ we have $T^{1+\e} \le M \ll T^{3/2}$,
and
$$
G_3 = \sup_{M\le x \le M',1\le t \le M^{1+\e}/T}\;\bigl|\sum_{h\le t}
\D _3(x,h)\bigr|.
$$
Moreover $t \le x^{{1\over3}+\delta}$ is satisfied, since
$$
t \le {M^{1+\e}\over T} \le M^{{1\over3}+\delta}
$$
because $M^{{2\over3}+\e-\delta} \le T$ for $\e < \delta$ and sufficiently
large $T$. Thus if (4.6) holds we have
$$
G_3M^{-1} \ll \sup_{M\le x\le2M,1\le t\le M^{1+\e}/T}\;tx^{{2\over3}+\e}M^{-1}
\ll M^{1+\e}T^{-1}M^{{2\over3}+\e}M^{-1} \ll T^\e,
$$
namely for $k = 3$ we have $a_3 = 0, b_3 = 1$ in (4.4).
Hence we obtain from (4.5)
$$
\int_0^T|\zt|^6\,\d t \;\ll_\e\;T^{1+\e},
$$
which was already shown to hold in [Iv4] if (5.10) is assumed. We also  have
from (4.5)

\medskip
{\bf Corollary 2}. {\it For a fixed integer $k \ge 3$ we have}
$$
\int_0^T|\zt|^{2k}\,\d t \;\ll_{k,\e}\;T^{1+\e} \leqno(4.7)
$$
{\it provided that (4.4) holds with $a_k = 0, b_k = 1$.}

\medskip
The bounds (4.5)  and (4.7) provide the means to bound the $2k$-th
moment of $|\zt|$. It remains to be seen, of course, whether there is
any hope of proving the $2k$--th    moment for $k \ge 4$ by (4.7),
namely whether $a_k = 0, b_k = 1$ can hold at all for sufficiently
large $k$. It would be highly interesting if one could even get
any non-trivial results concerning $a_k,b_k$ and improve unconditionally
the existing bounds (see Chapters 7 and 8 of [Iv1]) for the moments of $|\zt|$.
For more connections between divisor problems and power moments of $|\zt|$,
see e.g. the author's papers [Iv5] and [Iv6]. 

\medskip
\head
5. New  bounds for the sums  $\D _k(x,h)$ over the shift parameter $h$
\endhead
\medskip

S. Baier, T.D. Browning, G. Marasingha and L. Zhao  [BBMZ] recently proved that
$$
\sum_{h\le H} \Delta_3(N; h)
\ll_\e
N^{\varepsilon}\big(H^2+H^{1/2}N^{13/12}\big)
\qquad(1\le H\le N),
\leqno(5.1)
$$
$$
\Delta_3(N; h) = \sum_{N<n\leqslant2N}d_3(n)d_3(n+h) - NP_4(\log N;h),
$$
and if $N^{1/3+\varepsilon}\le H\le N^{1-\varepsilon}$, then
there exists $\delta = \delta(\e)>0$ for which
$$
\sum\limits_{h\le H} \left| \Delta_3(N; h) \right|^2
\;\ll_\e\; HN^{2-\delta(\e)}.
$$
Note that (5.1), in the interval $N^{1/6+\e}\le H \le N^{1-\e}$, gives
 an asymptotic formula for the averaged sum $\sum_{h\le H}D_3(N,h)$.

\medskip
Jie Wu and the author [IvWu] proved the following:

\medskip
THEOREM 2. {\it For a fixed integer $k\ge 3$ we have
$$
\sum_{h\le H}\D_k(N; h)
\ll_\e N^\e \big(H^2 + N^{1+\b_k}\big)
\qquad(1\le H\le N),
\leqno(5.2)
$$
where $\b_k$  is defined by
$$
\b_k
 := \inf\bigg\{\,b_k\;:\;\int_1^X |\D_k(x)|^2\d x \ll  X^{1+2b_k}\bigg\}
$$
and $\D_k(x)$ is the remainder term in the asymptotic formula for $\sum_{n\le x}d_k(n)$.}

\medskip
We have (see e.g., Chapter 13 of [Iv1])
$$
\sum_{n\le x}d_k(n) = xp_{k-1}(\log x) + \D_k(x),
$$
where
$$
p_{k-1}(\log x)
= \mathop{\roman {Res}}_{s=1}\bigg(\zeta^k(s) \frac{x^{s-1}}{s}\bigg),
$$
so that $p_{k-1}(z)$  is a polynomial of degree $k-1$ in $z$, all of whose coefficients
depend on $k$. In particular,
$$
p_1(z) \;=\;  z + 2\gamma -1\qquad(\gamma = -\G'(1)).
$$
It is known that $\b_k= (k-1)/(2k)$ for $k=2,3,4$, $\b_5 \le 9/20, \b_6 \le 1/2$, etc.
and $\b_k\ge (k-1)/(2k)$ for every $k\in \NN$.
It is conjectured that
$ \b_k= (k-1)/(2k)$  for every $ k\in \NN$, and this is equivalent to the
Lindel\"of Hypothesis  (that $\zeta(\frac{1}{2}+it)\ll_\e (|t|+1)^\e$).
From the Theorem 1 we obtain, for $1 \le H\le N$,
$$\eqalign{
\sum_{h\le H}\D_3(N; h) &\ll_\e N^\e \big(H^2 + N^{4/3}\big),
\cr
\sum_{h\le H}\D_4(N; h) &\ll_\e N^\e \big(H^2 + N^{11/8}\big),
\cr
\sum_{h\le H}\D_5(N; h) &\ll_\e N^\e \big(H^2 + N^{29/20}\big),
\cr
\sum_{h\le H}\D_6(N; h) &\ll_\e N^\e \big(H^2 + N^{3/2}\big).
\cr
}
$$
Since it is known that $\b_k <1$ for any $k$, this means
that the bound in (5.2) improves on the trivial bound $HN^{1+\e}$ in the
range $N^{\b_k+\e} \le H\le N^{1-\e}$. Our result thus  supports
the assertion that $\D_k(N; h)$ is really the error term in the asymptotic
formula for $D_k(N,h)$, as given above. In the case when $k=3$, we have
 an improvement on the result of Baier et al. when $H\ge N^{1/2}$.

\medskip
The basic idea of proof is to start from
$$
\eqalign{
&\sum_{h\le H} \Delta_k(N, h)
 = \sum_{N<n\le 2N} d_k(n) \sum_{h\le H} d_k(n+h) -
\sum_{h\le H} \int_N^{2N} \mathfrak{S}_k(x, h) \d x
\cr
& = M_k(N, H) + R_k(N, H) - \sum_{h\le H} \int_N^{2N} \mathfrak{S}_k(x, h) \d x,\cr}
$$
say. Here
$$
\mathfrak{S}_k(x,h) := \sum_{q=1}^{\infty} \frac{c_q(h)}{q^2} Q_k(x,q)^2,
$$
where $\mu(n)$ is the M\"obius function,
$c_q(h):=\sum_{d\mid (h,q)} d\mu(q/d)$ is the Ramanujan sum and
$Q_k(x,q)$ is a polynomial 
of degree $2k-2$
whose coefficients depend on $q$, and may be explicitly evaluated (see e.g., [BBMZ]).
Further we set
$$
\eqalign{
&M_k(N, H)
 := \sum_{N<n\le 2N} d_k(n) \mathop{\roman {Res}}_{s=1}
 \bigg(\zeta(s)^k \frac{(n+H)^s-n^s}{s}\bigg),
\cr&
R_k(N, H)
:= \sum_{N<n\le 2N} d_k(n) \big(\Delta_k(n+H)-\Delta_k(n)\big),\cr}
$$
and use complex integration to estimate $M_k(N,h)$ and then connect $R_k(N, H)$
to mean square estimates for $\D_k(x)$.
We have
$$
\eqalign{
M_k(N,H)
& = H\int_N^{2N}\big(\mathop{\roman {Res}}_{s=1}
\zeta(s)^kx^{s-1}\big)^2\d x
\cr
& \quad
+ O_\varepsilon\big(H^2N^\e + NH^{\alpha_k+\e} + N^{1+\beta_k+\e}\big)
\cr}
$$
and
$$
\sum\limits_{h\le H} \int_{N}^{2N} \mathfrak{S}_k(x,h)\d x
= H \int_{N}^{2N} \big(\mathop{\roman {Res}}
\limits_{s=1} \zeta(s)^kx^{s-1}\big)^2 \d x
+O_\varepsilon\big(N^{1+\varepsilon}\big).
$$
The constants $\alpha_k, \beta_k$ are defined as
$$
\alpha_k = \inf\Bigl\{\;a_k\;:\; \D_k(x) \ll a_k\;\Bigr\}
$$
and
$$
\b_k = \inf\Bigl\{\;b_k\;:\; \int_1^X\D_k^2(x)\d x \ll X^{1+2b_k}\;\Bigr\},
$$
and $(k-1)/(2k) \leqslant \b_k \leqslant \a_k < 1$ for $k = 2, \ldots\,.$
By completing the estimations one obtains the assertion of the theorem.

\medskip
\head
6. New results involving $\D(x+U)-\D(x)$
\endhead

The final topic will be a discussion (see (3.2)) of
$$
\D(x+U)-\D(x) = \sum_{x<n\le x+U}d(n) + O(Ux^\e)\qquad(1 \ll U \le x),\leqno(6.1)
$$
so that we are considering the divisor function $d(n)$ in ``short intervals''
$[x, x+U]$ if $U = o(x)$ as $x\to\infty$.
The interest in this topic comes from the work of M. Jutila
[Ju1], who  proved
that
$$
\eqalign{& \int\limits_T^{T+H}\Bigl(\D(x+U)-\D(x)\Bigr)^2\d x \cr&=
{1\over4\pi^2}\sum_{n\le {T\over2U}}{d^2(n)\over
n^{3/2}}\int\limits_T^{T+H} x^{1/2}\left|\exp\left(2\pi iU\sqrt{{n\over
x}}\,\right)-1\right|^2\d x + O_\e(T^{1+\e} +
HU^{1/2}T^\e),\cr}\leqno(6.2)
 $$
 for $1 \le U \ll T^{1/2} \ll H \le T$.
 From (6.2) one deduces ($a\asymp b$ means  $a\ll b\ll a$)
$$
\int_T^{T+H}\Bigl(\D(x+U)-\D(x)\Bigr)^2\d x
\asymp HU\log^3\left({\sqrt{T}\over U}\right)
\leqno(6.3)
$$
for $HU \gg T^{1+\e}$ and $T^\e \ll U \le \hf\sqrt{T}$. In [Ju2]
Jutila proved that the integral in (6.3) is
$$
\ll_\e T^\e(HU + T^{2/3}U^{4/3})\qquad(1 \ll H,U \ll X).
$$
In the case when $H=T$ the author [Iv7] improved (6.2) and proved

\medskip
THEOREM 3. {\it For $1 \ll U = U(T) \le \hf {\sqrt{T}}$ we have} ($c_3 = 8\pi^{-2}$)
$$\eqalign{
\int_T^{2T}\Bigl(\D(x+U)-\D(x)\Bigr)^2\d x & = TU\sum_{j=0}^3c_j\log^j
\Bigl({\sqrt{T}\over U}\Bigr) \cr&
+ O_\e(T^{1/2+\e}U^2) + O_\e(T^{1+\e}U^{1/2}).\cr}\leqno(6.4)
$$

\medskip
In (6.4) all the constants $c_j$  may be made explicit.
Note that, for $T^\e \le U = U(T) \le T^{1/2-\e}$ (6.4) is a true asymptotic formula.
From (6.4) one can deduce that, for $1 \ll U \le \hf {\sqrt{T}}$, we have ($c_3 = 8\pi^{-2}$)
$$\eqalign{
\sum_{T\le n \le 2T}\Bigl(\D(n+U)-\D(n)\Bigr)^2 & = TU\sum_{j=0}^3c_j\log^j
\Bigl({\sqrt{T}\over U}\Bigr) \cr&
+ O_\e(T^{1/2+\e}U^2) + O_\e(T^{1+\e}U^{1/2}).\cr}\leqno(6.5)
$$
The asymptotic formula (6.5) is a considerable improvement over a result of Coppola--Salerno [CoSa],
who had shown that ($T^\e \le U \le \hf \sqrt{T},\,L = \log T$)
$$
\sum_{T\le n \le 2T}\Bigl(\D(n+U)-\D(n)\Bigr)^2=
{8\over\pi^2}TU\log^3\Bigl({\sqrt{T}\over U}\Bigr)
+ O(TUL^{5/2}).
$$
The starting point for the above results is the explicit expression (see e.g., Chapter 3 of [Iv1])
$$
\D(x) = {1\over\pi\sqrt{2}}
x^{1\over4}\sum_{n\le N}d(n)n^{-{3\over4}}\cos(4\pi\sqrt{nx}
- {\txt{1\over4}}\pi) +
O_\e(x^{{1\over2}+\e}N^{-{1\over2}})\quad(2 \le N \ll x),
$$
which is flexible, since the parameter $N$ may be arbitrarily chosen. It is sometimes
called the {\it truncated Vorono{\"\i} formula}, in honour of G.F. Vorono{\"\i} [Vor],
who more than a century ago obtained an explicit formula for $\D(x)$ containing the
familiar Bessel functions.

\medskip
The most recent results on $\D(x+U)-\D(x)$ have been obtained by the author and W. Zhai
[IvZh].  We state just two of their theorems.

\medskip
THEOREM 4. {\it Suppose $\log^2 T\ll U\leqslant T^{1/2}/2,
T^{1/2}\ll H\leqslant T,$ then we have
$$
\eqalign{
 &\int_T^{T+H} \max_{0\leqslant u\leqslant U}\Bigl|\Delta(x+u)-\Delta(x)\Bigr|^2\d x \ll
HU{\Cal L}^5+T{\Cal L}^4\log{\Cal L}\cr
&   +H^{1/3}T^{2/3}U^{2/3}{\Cal L}^{10/3}(\log {\Cal L})^{2/3},\cr}
$$
where ${\Cal L} :=\log T.$}

This generalizes and sharpens a result of D.R. Heath-Brown \& K.-M. Tsang (1994) [HBTs].
From Theorem 3 we obtain then

\medskip
THEOREM 5. {\it
Suppose $T, U, H$ are large parameters and $C>1$ is
a large constant such that
$$T^{131/416+\varepsilon}\ll U\leqslant C^{-1}T^{1/2}{\Cal L}^{-5},
 \quad CT^{1/4}U{\Cal L}^5\log {\Cal L}\leqslant H\le T.$$
Then in the interval $[T, T+H]$ there are $\gg HU^{-1}$ subintervals
of length   $\gg U$ such that on each subinterval one has $
\pm\Delta(x) \ge c_{\pm} T^{1/4} $ for some $c_{\pm}>0.$}

\vfill
\eject
%\topglue1cm
\bigskip
\Refs
\bigskip

\item{[BBMZ]} S. Baier, T.D.  Browning, G. Marasingha and L. Zhao,
Averages of shifted convolutions of $d_3(n)$,
Proc. Edinb. Math. Soc. II. Ser. {\bf55}, No. 3,(2012) 551-576.

\item{[BeSh]} R. Bellman and H.N. Shapiro, On a problem in additive
number theory, Annals of Math. {\bf49}(1948), 333-340.

\item{[BEFSP]}
Y. Buttkewitz, C. Elsholtz, Christian; K. Ford and J.-C.
Schage-Puchta,
A problem of Ramanujan, Erd\H os and K\'atai on the iterated divisor function,
Int. Math. Res. Not. 2012, No. {\bf17}, 4051-4061(2012)

\item{[CoSa]} G. Coppola and S. Salerno, On the symmetry of the divisor
function in almost all short intervals, Acta Arith. {\bf113}(2004),
189-201.

\item{[EGIP]} P. Erd\H os, S.W. Graham,  A. Ivi\'c and C. Pomerance,
       On the number of divisors of $n!$,
Analytic Number Theory: Proceedings of a
       Conference in Honor of Heini Halberstam (Urbana, May 1995) Volume 1,
(eds. B.C. Berndt et al.), Birkh\"auser, Boston etc., 1996, 337-355.

\item{[ENS]} P. Erd\H os, J.-L. Nicolas and A. S\'ark\"ozy,
On large values of the divisor function,
Ramanujan J. {\bf2}, No.1-2, 225-245 (1998).

\item{[Erd1]} P. Erd\H os, Asymptotische Untersuchungen \"uber
die Anzahl der Teiler von $n$, Math. Annalen {\bf169}(1967), 230-238.

\item{[Erd2]} P. Erd\H os, On the sum $\sum_{n=1}^x d[d(n)]$,
Math. Stud. {\bf36}(1968), 227-229 (1969).

\item{[Erd3]} P. Erd\H os, ``Ramanujan and I'', Number Theory, Madras 1987, 81-92,
LNS in Mathematics 1395, Springer Verlag, Berlin etc., 1989.

\item{[ErHa]} P. Erd\H os and R.R. Hall, 
Values of the divisor function on short intervals,
J. Number Theory {\bf12}(1980), 176-187.

\item{[ErIv1]} P. Erd\H os and A. Ivi\'c, On the iterates of the enumerating
     function of finite Abelian groups,  Bulletin XCIX Acad. Serbe
    1989 Sciences Math\'ematiques
     No. {\bf17}, 13-22.

 \item{[ErIv2]} P. Erd\H os and A. Ivi\'c, The distribution of certain arithmetical functions
    at consecutive integers, Proceedings Budapest Conference in Number Theory
    July 1987, Coll. Math. Soc. J. Bolyai {\bf 51}, North-Holland, Amsterdam 1989,
    45-91.

\item{[ErKa1]} P. Erd\H os and I. K\'atai, On the growth of $d_k(n)$,
The Fibonacci Quarterly {\bf7}(1969), 267-274.

\item{[ErKa2]} P. Erd\H os and I. K\'atai, On the sum $\sum d_4(n)$, Acta. Scientarium Math.
(Szeged) {\bf30}(1969), 313-324.

\item{[ErMi]} P. Erd\H os and L. Mirsky,
The distribution of values of the divisor function $d(n)$,
Proc. Lond. Math. Soc. (3) {\bf2}(1952), 257-271.

\item{[Hal]} R.R. Hall,
The maximum value of the divisor function on short intervals, 
Bull. London Math. Soc. {\bf13}(1981), no. 2, 121–124.

\item{[HaTe]} R.R. Hall and G. Tenenbaum, On the local behaviour of some arithmetical functions,
Acta Arith. {\bf43}(1984), no. 4, 375-390.

\item{[HBTs]} D. R. Heath-Brown and K.Tsang, Sign changes of $E(t), \D(x)$ and $P(x)$,
J. of Number Theory {\bf49}(1994), 73-83.

\item{[Hux]} M.N. Huxley,
Exponential sums and lattice points III,
Proc. London Math. Soc., (3) {\bf 87}(2003), 591-609.

\item{[Ing]} A.E. Ingham, Some asymptotic formulae in the theory of numbers,
J. London Math. Soc. {\bf2}(1927), 205-208.

\item{[Iv1]} A. Ivi\'c, The Riemann zeta-function, John Wiley \&
Sons, New York, 1985 (2nd ed. Dover, Mineola, New York, 2003).

\item{[Iv2]} A. Ivi\'c, An asymptotic formula involving the enumerating function of finite
    Abelian groups, Publikacije Elektroteh. Fak. (Beograd) Ser. Mat. {\bf3}
   (1992), 61-66.

\item{[Iv3]} A. Ivi\'c:
On the ternary additive divisor problem and the sixth moment of the
zeta-function, in
 Sieve Methods, Exponential sums, and their Applications in
Number Theory (G.R.H. Greaves, G. Harman and M.N. Huxley, eds.),
Cambridge University Press, Cambridge, 1996, 205-243.

\item{[Iv4]} A. Ivi\'c, The general additive divisor problem and moments
of the zeta-function, in
``New Trends in Probability and Statistics'' Vol. 4, Analytic and
Probabilsitic Methods in Number Theory, eds. E. Laurin\v cikas etc.,
TEV Vilnius (Lithuania) and VSP (Utrecht \& Tokyo), 1997, pp. 69-89

\item{[Iv5]} A. Ivi\'c, On the Riemann zeta-function and the divisor problem,
Central European J. Math. {\bf(2)(4)} (2004), 1-15,  II, ibid.
{\bf(3)(2)} (2005), 203-214,  III, Annales Univ.
Sci. Budapest, Sect. Comp. {\bf29}(2008), 3-23,
and IV, Uniform Distribution Theory {\bf1}(2006), 125-135.

\item{[Iv6]} A. Ivi\'c, On the mean square of the zeta-function and
the divisor problem, Annales  Acad. Scien. Fennicae Mathematica {\bf23}(2007), 1-9.

\item{[Iv7]} A. Ivi\'c, On the divisor function and the Riemann zeta-function in short intervals,
The Ramanujan Journal: Volume 19, Issue 2 (2009), 207-224.

\item{[IvMo1]} A. Ivi\'c and Y. Motohashi, On the fourth power moment of the
    Riemann zeta-function,
    J. Number Theory {\bf51} (1995), 16-45.

\item{[IvMo2]} A. Ivi\'c and Y. Motohashi, On some estimates involving the binary
additive divisor problem, Quart. J. Math. Oxford (2) {\bf46}(1995), 471-483.

\item{[IvWu]} A. Ivi\'c and Jie Wu, On the General Additive Divisor Problem,
Proceedings of the Steklov Institute of Mathematics {\bf276}(2012), 140-148.

\item{[IvZh]} A. Ivi\'c and W. Zhai, On the Dirichlet divisor problem in short intervals,
The Ramanujan Journal, DOI 10.1007/s11139-013-9470-6 (2013).

\item{[Ju1]} M. Jutila, On the divisor problem for short intervals,
Ann. Univer. Turkuensis Ser. {\bf A}I {\bf186}(1984), 23-30.

\item{[Ju2]} M. Jutila, Riemann's zeta-function and the divisor problem,
Arkiv Mat. {\bf21}(1983), 75-96 and II, ibid. {\bf31}(1993), 61-70.

\item{[Ka1]} I. K\'atai, On the sum $\sum dd(f(n))$, Acta Sci. Math. {\bf29}(1968), 199-205.

\item{[Ka2]} I. K\'atai, On the iteration of the divisor function,
Publ. Math. debrecen {\bf16}(1969), 3-15 (1970).

\item{[Ka3]} I. K\'atai, On a problem of A. Ivi\'c, Mathematica Pannonica {\bf18}(2007),
11-18.

\item{[Kac]} J. Kaczorowski, Axiomatic theory of $L$-functions: the Selberg class, in
``Analytic Number Theory'' eds. A. Perelli and C. Viola, Springer Verlag,
Berlin-Heidelberg, 2006, pp. 133-209.

\item{[Meu]} T. Meurman, The binary  additive divisor problem, Number Theory Turku Proc., Turku
Symposium 199 (M. Jutila et al. editors), Walter de Gruyter, Berlin, 2001, 223-246.

\item{[Mot]} The binary  additive divisor problem, Annales sci. \'Ecole Norm. Sup., 4-\`eme
s\'erie, {\bf27}(1994), 529-572.

\item{[Ram]} S. Ramanujan, Highly composite numbers, Proc. London Math. Soc. 2 {\bf14}(1915), 347-409.

\item{[Ric]} H.-E. Richert, Versch\"arfung der Absch\"atzung beim Dirichletschen Teilerproblem,
Math. Zeit-\break schrift {\bf58}(1953), 204-218.

\item{[Sm1]} H. Smati, Sur un probl\`eme de S. Ramanujan, C. R., Math., Acad. Sci.
Paris {\bf340}, No. 1, 1-4, 2005.

\item{[Sm2]} H. Smati, Sur un probl\`eme d'Erd\H os et K\'atai,
Ann. Univ. Sci. Budap. Rolando E\"otv\"os, Sect. Comput. {\bf29}(2008), 213-238.

\item{[Vor]} G.F. Vorono{\"\i}, Sur une fonction transcendante et ses applications
\`a la sommation de quelques s\'eries, Ann. \'Ecole Normale (3){\bf21}(1904),
2-7-267 and ibid. 459-533.

\endRefs
%\vskip1cm

\enddocument

\bye